\documentstyle[12pt]{amsart}

\newcommand{\tor}{U_q({\frak g}_{tor})}
\newcommand{\torr}{U_q^{\phi'}({\frak g}_{tor})}
\newcommand{\affa}{U_q^{(1)}(\widehat{\frak {sl}}_{n+1})}
\newcommand{\affaa}{U_q^{(1)'}(\widehat{\frak {sl}}_{n+1})}
\newcommand{\affb}{U_q^{(2)}(\widehat{\frak {sl}}_{n+1})}
\newcommand{\affbb}{U_q^{(2)'}(\widehat{\frak {sl}}_{n+1})}

\theoremstyle{plain}
\newtheorem{thm}{Theorem}[subsection]
\newtheorem{prop}[thm]{Proposition}
\newtheorem{lemma}[thm]{Lemma}
\newtheorem{cor}[thm]{Corollary}
\newtheorem{defn}[thm]{Definition}
\theoremstyle{definition}

\newtheorem{rem}[thm]{Remark}

\newtheorem{ack}{Acknowledgment}

\numberwithin{equation}{subsection}

\begin{document}
\title[Quantum toroidal algebras]
 {Quantum toroidal algebras and their vertex representations}
\author{Yoshihisa Saito}
\address[ ]{Research Institute for Mathematical Sciences,
Kyoto University, Japan}
\thanks{The author is supported by the JSPS Research Fellowships for
Young Scientists.}
\maketitle
\begin{abstract}
 We construct the vertex representations of the quantum toroidal algebras
 $U_q({\frak {sl}}_{n+1,tor})$. In the 
 classical case the vertex representations are not irreducible. However in 
 the quantum
 case they are irreducible.\\ 
 {}For n=1, we construct
 a set of finitely many generators of $U_q({\frak {sl}}_{2,tor})$. 
\end{abstract}
\section{Introduction}
 The classical toroidal algebras have been studied by many authors [MEY], [S], 
 [Y], etc. Here ``classical'' means $q=1$.    
 The definition of the quantum toroidal algebras is given in [GKV]. They gave 
 a geometric realization of the quantum toroidal algebras without 
 any results on their representation theory. Recently
 Varagnolo and Vasserot [VV] proved a Schur-type duality between 
 representations of the quantum toroidal algebras and the double
 affine Hecke algebra introduced by Cherednik [C]. This is an analogue of 
 the duality between the quantum affine algebras and the affine Hecke 
 algebras given by Chari and Pressly [CP]. In [VV] only the representations
 of ``trivial central charge'' was studied. It is known that there are
 two subalgebras $\affa$ and $\affb$ of $U_q({\frak {sl}}_{n+1,tor})$ which 
 are isomorphic to $U_q(\widehat{\frak {sl}}_{n+1})$. We say the 
 $U_q({\frak {sl}}_{n+1,tor})$-module $M$ has the trivial central charge if $M$
 has level $0$ as both $\affa$-module and $\affb$-module. It is an 
 analogue of level 
 $0$ representations of the affine quantum algebras. In this paper we say 
 that $M$ has a level $(0,0)$ instead of the trivial central charge. The 
 first $0$ means that $M$ has a level $0$ as a $\affa$-module and the 
 second $0$ means that 
 $M$ has a level $0$ as a $\affb$-module.\\

 In this paper we try to consider an analogue of ``integrable 
 representations'' in the toroidal case. Let us recall the integrability of 
 quantum Kac-Moody modules. Let $U_q({\frak g})$ 
 be a quantum Kac-Moody
 algebra and $V$ a $U_q({\frak g})$-module. We say $V$ is integrable 
 if $V$ has a weight space decomposition and locally 
 nilpotent actions of 
 the Chevalley generators of $U_q({\frak g})$. Therefore the definition 
 of ``integrability'' needs 
 Chevalley-type generators. The toroidal algebras is defined through Drinfeld
 type of generators and its Chevalley-type generators are not known. Therefore
 we are not able to define the integrability at this moment. However, in the 
 affine case, Frenkel-Jing [FJ] realized the integrable
 representations with level 1 by the vertex representations. 
 Thus if there are ``vertex representations'' of quantum toroidal 
 algebras, they must be interesting examples of the integrable 
 representations still not defined.
 In the $q=1$ case, vertex representations of
 the toroidal algebras have been already considered by 
 Moody-Eswara Rao-Yokonuma [MEY].
 In this paper for we construct the $q$-analogue 
 of the representations
 defined by them for ${\frak g}={\frak {sl}}_{n+1}$ with level $(1,0)$ and 
 $(1,1)$. Therefore we give
 a new class of the representations of the quantum toroidal algebras. 
 In the $q=1$ case, the Fock modules
 are not irreducible over the Heisenberg algebra and the vertex 
 representations are not irreducible over the toroidal algebra.
 In the quantum case, it is not the case: the Fock modules are irreducible 
 and also the vertex 
 representations with level $(1,0)$ are irreducible. \\

 The algebra $\tor$ has infinitely many generators satisfying infinitely many
 relations (See \S 2). It is preferable that $\tor$ is written by 
 finitely many generators with finitely many relations. According to [GKV] 
 there are finitely many generators of
 $U_q({\frak {sl}}_{2,tor})$ but the relations among these generators are
 highly non-trivial. In this paper we give an explicit form of finitely many 
 generators of $U_q({\frak {sl}}_{2,tor})$ and closed relations of them 
 (See \S 4). They coincide with the generators by Vasserot 
 [V]. \\

\section{Definition of quantum toroidal algebras}
\subsection{Notations}
Let $\frak g$ be a complex semisimple Lie algebra of type $A_n$ and 
$\widehat{\frak g}$ an affine Kac-Moody Lie algebra of type $A_n^{(1)}$.
We denote their Cartan subalgebras by $\frak h$ and $\widehat{\frak h}$ 
respectively .
We denote by $\overline{\alpha}_1,\cdots,\overline{\alpha}_n$ the simple roots
of $\frak g$, by $\overline{h}_1,\cdots,\overline{h}_n$ the simple coroots of
$\frak g$, by $\overline{\Lambda}_1,\cdots,\overline{\Lambda}_n$ the 
fundamental
weights of $\frak g$,
by $\alpha_0,\cdots,\alpha_n$ the simple roots of $\widehat{\frak g}$,
by $h_0,\cdots,h_n$ the simple coroots of $\widehat{\frak g}$ and 
$\Lambda_0,\cdots,\Lambda_n$ the fundamental weights of $\widehat{\frak g}$.
Let $\overline{Q}=\oplus_{i=1}^n{\Bbb Z}\overline{\alpha}_i$ be the root
lattice of $\frak g$, $\overline{P}=\oplus_{i=1}^n{\Bbb Z}
\overline{\Lambda}_i$ the weight lattice of $\frak g$,
$Q=\oplus_{i=0}^n{\Bbb Z}\alpha_i$ the root lattice of $\widehat{\frak g}$ and
$P=\oplus_{i=0}^n{\Bbb Z}\Lambda_i\oplus{\Bbb Z}\delta$ the weight lattice
of $\widehat{\frak g}$. Here $\delta$ is the null root.
  
We denote the pairing of $\frak h$ and
${\frak h}^*$ ({\it resp}. $\widehat{\frak h}$ and $\widehat{\frak h}^*$) by 
$\langle\text{ },\text{ }\rangle$. The invariant bilinear form on $P$ is 
given by $(\alpha_i|\alpha_j)=-\delta_{ij-1}+2\delta_{ij}-\delta_{ij+1}$ and
$(\delta|\delta)=0$. The projection form $P$ to $\overline{P}$ is given by
$\overline{\Lambda_i}=\Lambda_i-\Lambda_0$ and $\overline{\delta}=0$.
\subsection{ }
We will give the definition of the quantum toroidal algebra $\tor$.
\begin{defn}
Let $M=(m_{ij})_{0\leq i,j\leq n}$ be a skew-symmetric $(n+1)\times
(n+1)$-matrix with integral coefficients and let 
$\kappa$ be an element of ${\Bbb Q}(q)^{\times}$.
$U_q({\frak g}_{tor})$ is an associated algebra over ${\Bbb Q}(q)$ with
generators :
$$E_{i,k},\text{ }F_{i,k},\text{ }H_{i,l},\text{ }K_i^{\pm},\text{ }
q^{\pm\frac12 c},q^{\pm d_1},q^{\pm d_2},$$
for $k\in {\Bbb Z}$, $l\in {\Bbb Z}\backslash \{0\}$ and $i=0,1,\cdots,n$.\\

We introduce $K^{\pm}_{i,k}$ as the Fourier components of the following
generating functions:
$$K_i^+(z)=\sum_{k\geq 0}K^+_{i,k}z^{-k}=K_i^+\exp((q-q^{-1})\sum_{k\geq 1}
H_{i,k}z^{-k}),$$
$$K_i^-(z)=\sum_{k\leq 0}K^-_{i,k}z^{-k}=K_i^-\exp(-(q-q^{-1})\sum_{k\geq 1}
H_{i,-k}z^{k}).$$
The defining relations of $\tor$ are then written as follows:
\begin{equation}
 q^{\pm\frac12 c}\text{ are central,}
\end{equation}
\begin{equation}
 K_i^+K_i^{-}=K_i^{-}K_i^+=1,
\end{equation}
\begin{equation}
 [K_i^{\pm},K_j^{\pm}]=0,
\end{equation}
\begin{equation}
 [K_i^{\pm},H_{j,l}]=0,
\end{equation}
\begin{equation}
 [H_{i,k},H_{j,l}]=\delta_{k+l,0}\frac{1}{k}[k\langle h_i,\alpha_j\rangle]
 \frac{q^{kc}-q^{-kc}}{q-q^{-1}}\kappa^{-km_{ij}},
\end{equation}
\begin{equation}
 [q^{\pm d_i},K_j^{\pm}]=0,
\end{equation}
\begin{equation}
 q^{d_1}H_{j,l}q^{-d_1}=q^lH_{j,l},
\end{equation}
\begin{equation}
 [q^{\pm d_2},H_{j,l}]=0,
\end{equation}
\begin{equation}
 q^{d_1}E_{j,k}q^{-d_1}=q^kE_{j,k},
\end{equation}
$$
 q^{d_1}F_{j,k}q^{-d_1}=q^kF_{j,k},
$$ 
\begin{equation}
 q^{d_2}E_{j,k}q^{-d_2}=q^{\delta_{j0}}E_{j,k},
\end{equation} 
$$
 q^{d_2}F_{j,k}q^{-d_2}=q^{-\delta_{j0}}F_{j,k},
$$
\begin{equation}
 K_i^+E_{j,k}K_i^{-}=q^{\langle h_i,\alpha_j\rangle}E_{j,k},
\end{equation}
\begin{equation*}
 K_i^+F_{j,k}K_i^{-}=q^{-\langle h_i,\alpha_j\rangle}F_{j,k},
\end{equation*}
\begin{equation}
 [H_{i,k},E_{j,l}]=\frac{1}{k}[k\langle h_i,\alpha_j\rangle]q^{-\frac12 |k|c}
 \kappa^{-km_{ij}}E_{j,k+l},
\end{equation}\label{fr}
\begin{equation*}
 [H_{i,k},F_{j,l}]=-\frac{1}{k}[k\langle h_i,\alpha_j\rangle]q^{\frac12 |k|c}
 \kappa^{-km_{ij}}F_{j,k+l},
\end{equation*}
\begin{equation}
 \kappa^{m_{ij}}E_{i,k+1}E_{j,l}-q^{\langle h_i,\alpha_j\rangle}
 \kappa^{m_{ij}}E_{j,l}E_{i,k+1}
 =q^{\langle h_i,\alpha_j\rangle}
 E_{i,k}E_{j,l+1}-E_{j,l+1}E_{i,k},
\end{equation}
\begin{equation*}
 \kappa^{m_{ij}}F_{i,k+1}F_{j,l}-q^{-\langle h_i,\alpha_j\rangle}
 \kappa^{m_{ij}}{}F_{j,l}F_{i,k+1}
 =q^{-\langle h_i,\alpha_j\rangle}
 {}F_{i,k}F_{j,l+1}-F_{j,l+1}F_{i,k},
\end{equation*}
\begin{equation}
 [E_{i,k},F_{j,l}]=\delta_{i,j}\frac{1}{q-q^{-1}}\{q^{\frac12 (k-l)c}
 K_{i,k+l}^+-q^{\frac12 (l-k)c}K_{i,k+l}^-\},
\end{equation} 
\begin{equation}
 \sum_{\sigma\in {\frak S}_m}\sum_{r=0}^m(-1)^r
 \begin{bmatrix}m\\r\end{bmatrix}
 E_{i,k_{\sigma(1)}}\cdots E_{i,k_{\sigma(r)}}E_{j,l}E_{i,k_{\sigma(r+1)}}
 \cdots E_{i,k_{\sigma(m)}}=0,
\end{equation}
\begin{equation*}
 \sum_{\sigma\in {\frak S}_m}\sum_{r=0}^m(-1)^r
 \begin{bmatrix}m\\r\end{bmatrix}
 {}F_{i,k_{\sigma(1)}}\cdots F_{i,k_{\sigma(r)}}F_{j,l}F_{i,k_{\sigma(r+1)}}
 \cdots F_{i,k_{\sigma(m)}}=0,
\end{equation*}
\qquad\qquad\qquad\qquad\qquad\qquad\qquad\qquad\qquad\qquad
\qquad\qquad\qquad\qquad for $i\ne j$,

where $m=1-\langle h_i,\alpha_j\rangle$.\\

In these relations we denote 
 $[k]=\frac{q^k-q^{-k}}{q-q^{-1}}$, $[n]!=\prod_{k=1}^n[k]$, 
 $\begin{bmatrix}m\\r\end{bmatrix}=\frac{[m]!}{[r]![m-r]!}$.
\end{defn}
\subsection{}
Let $U_q^{''}({\frak g}_{tor})$ be the subalgebra of $\tor$ generated by 
$E_{i,k},F_{i,k},K_i^{\pm},H_{i,l}$, $q^{\pm\frac12 c}$. Let 
$U_q^{\phi'}({\frak g}_{tor})$ ({\it resp} $U_q^{'\phi}({\frak g}_{tor})$) 
be the subalgebra
generated by $U_q^{''}({\frak g}_{tor})$ and $q^{\pm d_1}$ 
({\it resp} $q^{\pm d_2}$). 
Let $\affaa$ be the subalgebra generated
by $E_{i,k},F_{i,k},K_i^{\pm},H_{i,l}$, $q^{\pm\frac12 c}$ 
$(1\leq i \leq n,k\in{\Bbb Z},
l\in{\Bbb Z}\backslash\{0\})$ and $\affa$ the subalgebra generated by
$\affaa$ and $q^{\pm d_1}$. 
Let $\affbb$ the subalgebra generated by $E_{i,0},F_{i,0},K_i^{\pm}$
$(0\leq i \leq n)$ and $\affb$ the subalgebra generated by
$\affbb$ and $q^{\pm d_2}$.
By the definition $\affb$ is isomorphic to 
$U_q(\widehat{\frak {sl}}_{n+1})$ and $\affbb$ is isomorphic to 
$U_q'(\widehat{\frak {sl}}_{n+1})$.\\ 

The following are straight forward.
\begin{lemma}
 For $1\leq i \leq n$ let $\overline{E_{i,k}}=E_{i,k}\kappa^{\sum_{j=1}^i
 km_{j-1j}}$, $\overline{F_{i,k}}=F_{i,k}\kappa^{\sum_{j=1}^i
 km_{j-1j}}$, $\overline{H_{i,l}}=H_{i,l}\kappa^{\sum_{j=1}^i
 lm_{j-1j}}$. Then the relations between $\overline{E_{i,k}}$, 
 $\overline{F_{i,k}}$, $\overline{H_{i,l}}$ and $K_i^{\pm}$ are precisely 
 the relations of Drinfeld generators of $U_q(\widehat{\frak {sl}}_{n+1})$.
 That is, $\affa$ is isomorphic to $U_q(\widehat{\frak {sl}}_{n+1})$
 and $\affaa$ is isomorphic to $U_q'(\widehat{\frak {sl}}_{n+1})$.
\end{lemma}

\begin{lemma}
 Let $K_{\delta}^{\pm}=\prod_{i=0}^nK_i^{\pm}$. Then $K_{\delta}^{\pm}$ 
 are the central elements of $\tor$.
\end{lemma}

Note that $q^{\pm c}$ is the central elements of $\affa$ and 
$K_{\delta}^{\pm}$ the central elements of $\affb$.
\section{Vertex representations}
\subsection{Heisenberg algebras}
 In this section we shall give the vertex representations of $\tor$.
 We assume $c=1$.\\

 Consider a ${\Bbb Q}(q)$-algebra $S_n$ generated by $H_{i,l}$
 $(0\leq i \leq n,\text{ } l\in
 {\Bbb Z}\backslash \{0\})$ satisfying:
 \begin{equation}
  [H_{i,k},H_{j,l}]=\delta_{k+l,0}\frac1k [k\langle h_i,\alpha_j\rangle]
  \frac{q^{k}-q^{-k}}{q-q^{-1}}\kappa^{-km_{ij}}.
 \end{equation}

 We call $S_n$ the Heisenberg algebra.

 Let $S_n^+$ ({\it resp.} $S_n^-$) be the subalgebra of $S_n$ generated by
 $H_{i,l}$ $(0\leq i \leq n,\text{ } l>0)$ ({\it resp.} $0\leq i \leq n,
 \text{ }l<0$). 

 We introduce the Fock space 
 $${\cal F}_n=S_nv_0$$
 with the defining relations:
 \begin{equation}
  H_{i,l}v_0=0, \qquad\text{for }l>0,
 \end{equation}
 \begin{equation}
  q^{\frac12 c}v_0=q^{\frac12}v_0.
 \end{equation}
 Note that ${\cal F}_n$ is a free $S_n^-$-module of rank $1$.\\

 Let $\Bbb F$ be a field of characteristic zero and  let $\frak a$ be 
 an associative $\Bbb F$-algebra generated by $x_p$, $y_p$ ($p\in 
 {\Bbb Z}_{> 0}$), $z$ and its inverse $z^{-1}$ with the following relations:
 $$[x_p,z]=[y_p,z]=0,$$
 $$[x_p,x_r]=[y_p,y_r]=0,$$
 $$[x_p,y_r]=\delta_{pr}z.$$
 Let ${\frak a}^+$ ({\it resp.} ${\frak a}^-$) be the subalgebra of 
 $\frak a$ generated by $x_p$ ({\it resp.} $y_p$). We set ${\frak b}=
 {\frak a}^+\otimes {\Bbb F}[z,z^{-1}]$. This is a maximal abelian 
 subalgebra of
 $\frak a$. Fix a nonzero scalar $\lambda\in {\Bbb F}^{\times}$. Let 
 ${\Bbb F}_{\lambda}$ be the one-dimensional space $\Bbb F$ viewed
 as a ${\frak b}$-module by:
 $$z\cdot 1=\lambda,\quad {\frak a}^+\cdot 1=0.$$
 Let $F(\lambda)$ be the induced ${\frak a}$-module 
 $$F(\lambda)=\text{Ind}^{\frak a}_{\frak b}{\Bbb F}_{\lambda}={\frak a}
 \otimes_{\frak b}
 {\Bbb F}_{\lambda}.$$
 By the defining relations of ${\frak a}$ we obtain an $\Bbb F$-linear
 isomorphism
 $$F(\lambda)\cong {\frak a}^-.$$
 Since ${\frak a}^-$ is abelian we may regard it as the algebra of polynomials
 in the variables $y_1$, $y_2$, $\cdots$. Then we see 
 that $z$ acts on ${\frak a}^-\cong
 {\Bbb F}[y_1,y_2,\cdots]$ by the multiplication of $\lambda$, $x_p$ acts by
 $\lambda\frac{\partial}{\partial y_p}$. By this realization we immediately
 have the following lemma.
 \begin{lemma}
  $F(\lambda)$ is an irreducible $\frak a$-module.
 \end{lemma}

 {}Fix an skew-symmetric $(n+1)\times (n+1)$-matrix with integral coefficients 
 $M=(m_{ij})_{0\leq i,j\leq n}$. We say that $\kappa\in {\Bbb Q}(q)$ is 
 generic with respect to 
 $M$ if for any $k\in{\Bbb Z}_{>0}$ the matrix 
 $([k\langle h_i,\alpha_j\rangle]\kappa^{-km_{ij}})$ is invertible.

 Note that if $n=1$ any $\kappa$ is not generic with respect to any $M$. 
 Since the matrix $([k\langle h_i,\alpha_j\rangle])_{0\leq i,j \leq n}$ is
 invertible for $n>1$, there exists a generic $\kappa$ for $n>1$.  
 \begin{lemma}
  $(1)$ For a fixed $M$, we assume that $\kappa\in {\Bbb Q}(q)$ is 
  generic with respect to $M$ (in particular $n>1$). Then ${\cal F}_n$ is an 
  irreducible $S_n$-module.\\
  $(2)$ ${\cal F}_1$ is not irreducible.
 \end{lemma}

 \begin{pf}
 (1) Set $G(k)=(g(k)_{ij})=([k\langle h_i,\alpha_j\rangle]\kappa^{-km_{ij}})$. 
 Since $\kappa$ is generic with respect to $M$ there exists its inverse
 $G(k)^{-1}=(g(k)^{ij})$ for any $k$. Note that by the definition 
 $\sum_{0\leq s\leq n}g(k)^{is}g(k)_{sj}=\delta_{ij}$.

 We set 
 $$\tilde{H}_{i,k}=
 \begin{cases}
   \sum_{0\leq s\leq n}\frac{k}{[k]}g(k)^{is}H_{s,k}, 
     & \qquad \text{for $k>0$},\\
   H_{i,k}, & \qquad \text{for $k<0$}.
 \end{cases}$$
 Then we have
 $$[\tilde{H}_{i,k},\tilde{H}_{j,l}]=[\tilde{H}_{i,-k},\tilde{H}_{j,-l}]=0,$$
 $$[\tilde{H}_{i,k},\tilde{H}_{j,-l}]=\delta_{k,l}\delta_{ij},$$
 for $k,l>0$.
 Since all $G(k)$ are regular, $\tilde{H}_{i,k}$ $(0\leq i\leq n,\text{ }k>0)$
 generate $S_n^+$.\\ 

 We shall use Lemma 3.1.1. Put ${\Bbb F}={\Bbb Q}(q)$, $\frak a=S_n$, 
 ${\frak a}^{\pm}=S_n^{\pm}$, $\lambda=1$, $x_p=\tilde{H}_{i,k}$ $(k>0)$, 
 $y_r=\tilde{H}_{j,l}$ $(l<0)$ where $p=(k-1)(n+1)+i+1$ and 
 $r=(-l-1)(n+1)+j+1$.   
   
 Then it is clear that $F(1)={\cal F}_n$. By Lemma 3.1.1 we conclude that
 ${\cal F}_n$ is an irreducible $S_n$-module.
 
 (2) It is easy to see that $\kappa^{-km_{01}}H_{0,-k}+H_{1,-k}$ is a central
 element of $S_1$ for each $k\in{\Bbb Z}_{>0}$. Therefore ${\cal F}_1$ has
 infinitely many singular vectors.
 \end{pf}
\subsection{Construction of level (1,0) modules}
 We assume $n>1$ and $M$ is an $(n+1)\times (n+1)$-matrix defined as follows;
 $$
M=\begin{pmatrix}
      0 &     -a &      0 & \hdots &      0 &      a\\
      a &      0 &     -a & \hdots &      0 &      0\\
      0 &      a &      0 & \hdots &      0 &      0\\
 \vdots & \vdots & \vdots & \ddots & \vdots & \vdots\\
      0 &      0 &      0 & \hdots &      0 &     -a\\
     -a &      0 &      0 & \hdots &      a &      0
\end{pmatrix}
$$
here $a\in {\Bbb Z}$.\\

 Note that one can rewrite $\overline{P}=\oplus_{i=2}^n{\Bbb Z}
 \overline{\alpha_i}
 \oplus{\Bbb Z}\overline{\Lambda}_n$.
 We introduce a
 twisted version of the group algebra ${\Bbb Q}(q)[\overline{P}]$ by 
 ${\Bbb Z}/2{\Bbb Z}$. We denote it by
 ${\Bbb Q}(q)\{\overline{P}\}$. This
 is the ${\Bbb Q}(q)$-algebra generated by symbols $e^{\overline{\alpha_2}}$, 
 $e^{\overline{\alpha_3}},
 \cdots,\text e^{\overline{\alpha_n}}$, $e^{\overline{\Lambda_n}}$ 
 which satisfy the 
 following relations:
 \begin{equation}
  e^{\overline{\alpha_i}}e^{\overline{\alpha_j}}=(-1)^{\langle\overline{h_i},
 \overline{\alpha_j}\rangle}e^{\overline{\alpha_j}}e^{\overline{\alpha_i}}
 \end{equation}
 \begin{equation}
  e^{\overline{\alpha_i}}e^{\overline{\Lambda_n}}=(-1)^{\delta_{in}}
  e^{\overline{\Lambda_n}}e^{\overline{\alpha_i}}.
 \end{equation}
 For $\overline{\alpha}=\sum_{i=2}^n m_i\overline{\alpha_i}+m_{n+1}
 \overline{\Lambda_n}$ we denote
 $e^{\overline{\alpha}}=(e^{\overline{\alpha_2}})^{m_2}
 (e^{\overline{\alpha_3}})^{m_3}
 \cdots (e^{\overline{\alpha_n}})^{m_n}
 (e^{\overline{\Lambda_n}})^{m_{n+1}}$. For example 
 $e^{\overline{\alpha_1}}=e^{-2\overline{\alpha_2}}e^{-3\overline{\alpha_3}}
 \cdots e^{-n\overline{\alpha_n}}
 e^{(n+1)\overline{\Lambda_n}}$, $e^{\overline{\Lambda_i}}=
 e^{-\overline{\alpha_{i+1}}}
 e^{-2\overline{\alpha_{i+2}}}\cdots e^{-(n-i)\overline{\alpha_n}}
 e^{(n-i+1)\overline{\Lambda_n}}$
 where $\overline{\Lambda_i}$ is the $i$-th fundamental weight.
 We denote $\overline{\alpha}_0=-\sum_{i=1}^n
 \overline{\alpha}_i$ and $\overline{h}_0=-\sum_{i=1}^n
 \overline{h}_i$.\\

 Note that $\langle h_i,\alpha_j \rangle=\langle\overline{h}_i,
 \overline{\alpha}_j\rangle$ for $0\leq i,j\leq n$.\\

 We denote by ${\Bbb Q}(q)\{\overline{Q}\}$ the subalgebra of
 ${\Bbb Q}(q)\{\overline{P}\}$ generated by $e^{\overline{\alpha_i}}$ 
 ($1\leq i \leq n$). \\

 Set 
 $$W(p)_n={\cal F}_n\otimes {\Bbb Q}(q)\{\overline{Q}\}
 e^{\overline{\Lambda_p}}\qquad\text{for }0\leq p\leq n.$$
 We define the operators $H_{i,l}$, $e^{\overline{\alpha}}$ 
 $(\overline{\alpha}\in \overline{Q})$,
 $\partial_{\overline{\alpha}_i}$, $z^{H_{i,0}}$, $d$ on $W(p)_n$ for 
 $i=0,1,\cdots,n$
 as follows:\\
 for $v\otimes e^{\overline{\beta}}=H_{i_1,-k_1}\cdots H_{i_N,-k_N}v_0
 \otimes e^{\sum_{j=1}^n m_j\overline{\alpha}_j+\overline{\Lambda}_p}
 \in W(p)_n$,
 $$H_{i,l}(v\otimes e^{\overline{\beta}})=(H_{i,l}v)\otimes 
 e^{\overline{\beta}},$$
 $$e^{\overline{\alpha}}(v\otimes e^{\overline{\beta}})=v\otimes 
 e^{\overline{\alpha}}e^{\overline{\beta}},$$
 $${\partial_{\overline{\alpha_i}}}(v\otimes e^{\overline{\beta}})=
 {\langle\overline{h_i},\overline{\beta}\rangle}v\otimes 
 e^{\overline{\beta}},$$
 $$z^{H_{i,0}}(v\otimes e^{\overline{\beta}})
  =z^{\langle \overline{h_i},\overline{\beta}\rangle}
  \kappa^{\frac12 \sum_{j=1}^n\langle \overline{h_i},m_j\overline{\alpha}_j
  \rangle m_{ij}}v\otimes e^{\overline{\beta}},$$
  $$d(v\otimes e^{\overline{\beta}})=(-\sum_{s=1}^{N}k_s
 -\frac{(\overline{\beta}|\overline{\beta})}{2}
 +\frac{(\overline{\Lambda_p}|\overline{\Lambda_p})}{2})v\otimes 
 e^{\overline{\beta}}.$$

We have the following lemma.

\begin{lemma}
 As operators on $W(p)_n$,
 $$e^{\overline{\alpha}_i}e^{\overline{\alpha}_j}
 =(-1)^{\langle \overline{h}_i,\overline{\alpha}_j\rangle}
 e^{\overline{\alpha}_j}e^{\overline{\alpha}_i},$$
 $$e^{\overline{\alpha}_i}q^{\partial_{\overline{\alpha}_j}}
 =q^{\langle \overline{h}_i,\overline{\alpha}_j\rangle}
 q^{\partial_{\overline{\alpha}_j}}e^{\overline{\alpha}_i},$$
 $$e^{\overline{\alpha}_i}z^{\partial_{\overline{\alpha}_j}}
 =z^{\langle \overline{h}_i,\overline{\alpha}_j\rangle}
 z^{\partial_{\overline{\alpha}_j}}e^{\overline{\alpha}_i}$$
 for $0\leq i,j \leq n$.
\end{lemma}

 We introduce the following generating functions:
 $$E_i(z)=\sum_{k\in {\Bbb Z}}E_{i,k}z^{-k},$$
 $$F_i(z)=\sum_{k\in {\Bbb Z}}F_{i,k}z^{-k}.$$

\begin{prop}
 Assume $c=1$ and $n>1$. Then
 for each $p$ and $\kappa$, the following action gives a $\torr$-module structure
 on $W(p)_n$:
 $$q^{\frac12 c}\mapsto q^{\frac12},$$
 $$q^{d_1}\mapsto q^{d},$$
 $$E_i(z)\mapsto\exp(\sum_{k\geq 1}\frac{H_{i,-k}}{[k]}(q^{-1/2}z)^k)
         \exp(\sum_{k\geq 1}-\frac{H_{i,k}}{[k]}(q^{1/2}z)^{-k})
         e^{\overline{\alpha_i}}z^{H_{i,0}+1},$$
 $$F_i(z)\mapsto\exp(\sum_{k\geq 1}-\frac{H_{i,-k}}{[k]}(q^{1/2}z)^k)
         \exp(\sum_{k\geq 1}\frac{H_{i,k}}{[k]}(q^{-1/2}z)^{-k})
         e^{-\overline{\alpha_i}}z^{-H_{i,0}+1},$$
 $$K_i^+(z)\mapsto \exp((q-q^{-1})
   \sum_{k\geq 1}H_{i,k}z^{-k})q^{\partial_{\overline{\alpha_i}}},$$
 $$K_i^-(z)\mapsto \exp(-(q-q^{-1})
   \sum_{k\geq 1}H_{i,-k}z^{k})q^{-\partial_{\overline{\alpha_i}}}$$
 for $0\leq i\leq n$.
\end{prop}

The proof will be given in Appendix.\\

We have immediately the following lemma.
\begin{lemma}
 The $\torr$-module $W(p)_n$ is cyclic: 
 $W(p)_n=\torr (v_0\otimes e^{\overline{\Lambda}_p}).$
\end{lemma}

\begin{thm}
 If $n>1$ and $\kappa$ is generic then $W(p)_n$ is irreducible for any $p$.
\end{thm}

\begin{pf}
 Since ${\cal F}_n$ is irreducible with respect to the action of $S_n$, 
 it is enough to show that for any non-zero
 $v=v_0\otimes \sum_{\overline{\alpha}\in\overline{Q}}
 a_{\overline{\alpha}}
 e^{\overline{\alpha}}e^{\overline{\Lambda}_i}$ ($a_{\overline{\alpha}}
 \in {\Bbb Q}(q)$) there exists $X\in \torr$ 
 such that $Xv=v_0\otimes e^{\overline{\Lambda}_p}$. 
 Let $\overline{S_n}$ be the 
 subalgebra of $S_n$ generated by $H_{i,l}$ $(1\le i\le n,l\in{\Bbb Z}
 \backslash\{0\})$ and $q^{\frac12 c}$. Let $\overline{\cal F}_n$
 be the $\overline{S_n}$-submodule of ${\cal F}_n$ generated by $v_0\otimes 
 e^{\overline{\Lambda}_i}$,
 and let $\overline{W(p)_n}=\overline{\cal F}_n\otimes 
 {\Bbb Q}(q)\{\overline{Q}\}e^{\overline{\Lambda_p}}$. As already known
 $\overline{W(p)_n}$ is an irreducible $\affa$-module. It is obvious that
 $v\in \overline{W(p)_n}$. Therefore there exists $X\in \affa\subset\torr$
 such that $Xv=v_0\otimes e^{\overline{\Lambda}_p}$. 
\end{pf}

\begin{rem}
 Since $\affa$ and $\affbb$ are subalgebras of $\torr$ we can regard $W(p)_n$
 as a $\affa$-module or as a $\affbb$-module. As a $\affa$-module, $W(p)_n$ is
 a level $1$ module. On the other hand it is a level $0$ module as a 
 $\affbb$-module. 
\end{rem}
\subsection{On the structure of level (1,0) modules}
 In this subsection we will study the $\affa$-module structure of $W(p)_n$.\\

 Let $M$ be a $\affa$-module.
 We denote the character of $M$ by $\text{ch}_M$.\\

 Let $L(\Lambda_p)$ be the irreducible highest weight 
 $U_q(\widehat{\frak {sl}}_{n+1})$-module with highest weight $\Lambda_p$.
 Note that the following identity
 holds:
 \begin{equation}
 \text{ch}_{L(\Lambda_p)}=\frac{e^{\Lambda_p}\sum_{\alpha\in \overline{Q}}
 e^{\alpha-(\frac12(\alpha|\alpha)-(\alpha|\overline{\Lambda_p}))\delta}}
 {\varphi(e^{-\delta})^{n}}.
 \end{equation}
 Here $\varphi(x)=\prod_{k>0}(1-x^k)$.\\

 We denote $\delta_1$ by the null root of $\affa$. 

 By the definition of $W(p)_n$ and (3.3.1) it is immediate to see the 
 following proposition.

 \begin{prop}
  \begin{align*}
   \text{ch}_{W(p)_n} 
   & = \frac{e^{\Lambda_p}\sum_{\alpha\in \overline{Q}}
         e^{\alpha-(\frac12(\alpha|\alpha)-(\alpha|\overline{\Lambda_p}))
         \delta_1}}{\varphi(e^{-\delta})^{n}}\\
   & = \frac{\text{ch}_{L(\Lambda_p)}}{\varphi(e^{-\delta_1})}.
  \end{align*}
 \end{prop}

 \begin{lemma}
  {}For each $l\in {\Bbb Z}\backslash\{0\}$ there 
  exist $\widehat{H}_l=\sum_{i=0}^{n}a_{i,l}H_{i,l}$
  $(a_{i,l}\in {\Bbb Q}(q))$ such that
  $$
  [\widehat{H}_l,H_{j,k}]=0
  $$
  for any $1\leq j \leq n$ and $k\in {\Bbb Z}\backslash\{0\}$. 
  Moreover such $\widehat{H}_l$ is unique up to scalar. 
 \end{lemma}

 \begin{pf}
  Note that $\kappa=1$. The rank of $n\times(n+1)$-matrix 
  $([l\langle h_i,\alpha_j\rangle])_{1\leq i \leq n,0\leq j \leq n}$ is equal
  to $n$. The lemma follows form this fact immediately.
 \end{pf}

 By the definition of $\widehat{H}_l$ we have
 \begin{equation}
  [\widehat{H}_k,\widehat{H}_l]=\delta_{k+l,0}\gamma_k,
 \end{equation}
 where $\gamma_k\in {\Bbb Q}(q)$.
 We fix a normalization of $\widehat{H}_l$ by putting $\gamma_k=1$ for all $k$.

 Let $\widehat{S}_n$ be the subalgebra of $S_n$ generated by 
 $\widehat{H}_l$. By the definition, $\widehat{S}_n$ acts on $W(p)_n$.\\ 

 The following two lemmas are easy to see.

 \begin{lemma}
  {}For $l>0$, $\widehat{H}_l(v_0\otimes e^{\overline{\beta}})=0$.
 \end{lemma}

 \begin{lemma}
  The action of $\affaa$ on $W(p)_n$ commutes with the action of 
  $\widehat{S}_n$.
 \end{lemma}

 Let $\widehat{S}_n^-$ be the subalgebra of $\widehat{S}_n$ generated by 
 $\widehat{H}_l$ $(l<0)$. 

 \begin{prop}
  As $\affa$-module
  $$W(p)_n\cong L(\Lambda_p)^{\oplus\infty}.$$
 \end{prop}

 \begin{pf}
  Set $\text{deg}(\widehat{H}_{k})=k$. Let $M_k=M_k(\widehat{H}_{-1},
  \widehat{H}_{-2},\cdots)$ be a monomial of degree $k$ in variables
  $\widehat{H}_{-1},\widehat{H}_{-2},\cdots$. 
  Then by the above two lemmas, $M_kv_0\otimes e^{\overline{\Lambda_p}}$
  is a singular vector of $\affa$-module $W(p)_n$. Let $W_{M_k}$ be the
  $\affa$-submodule which is generated by $M_kv_0\otimes 
  e^{\overline{\Lambda_p}}$. Then by the definition of the action of $\affa$
  on $W(p)_n$, we have
  $$W_{M_k}\cong L(\Lambda_p-k\delta)\cong L(\Lambda_p).$$  

  The vectors $\{M_kv_0\otimes e^{\overline{\Lambda_p}}\}$ are linearly 
  independent. The number of the monomials of degree $k$ is equal to the 
  $k$-th partition number $p(k)$. Therefore there is a $\affa$-submodule $W$ 
  of $W(p)_n$ which is isomorphic to $\oplus_{k\geq 0}
  L(\Lambda_p-k\delta)^{\oplus p(k)}$. By Proposition 3.3.1 it coincides with
  $W(p)_n$. This completes proof.   
 \end{pf}  

 By Lemma 3.3.4 and the proof of Proposition 3.3.5, the following corollary 
 follows 
 immediately.
 \begin{cor}
  As $\affaa\otimes \widehat{S}_n$-module $W(p)_n$ is isomorphic to
  $L(\Lambda_p)\otimes\widehat{S}_n^-$.
 \end{cor} 
\subsection{Construction of level (1,1) modules}
 We introduce a twisted version of the group algebra ${\Bbb Q}(q)[Q]$ by
 ${\Bbb Z}\backslash 2{\Bbb Z}$. We denoted it by ${\Bbb Q}(q)\{Q\}$.
 This is the ${\Bbb Q}(q)$-algebra generated by symbols
 $e^{\alpha_0}$,$e^{\alpha_1},\cdots,e^{\alpha_n}$ which satisfy the 
 following relations:
 \begin{equation}
  e^{\alpha_i}e^{\alpha_j}=(-1)^{\langle h_i,\alpha_j\rangle}
  e^{\alpha_j}e^{\alpha_i}.
 \end{equation}
 Similarly to \S 3.2, we denote $e^{\alpha}=(e^{\alpha_0})^{m_0}
 (e^{\alpha_1})^{m_1}\cdots
 (e^{\alpha_n})^{m_n}$ for $\alpha=\sum_{i=0}^nm_i\alpha_i\in Q$.

 Let 
 $$V(p)_n={\cal F}_n\otimes {\Bbb Q}(q)\{Q\}e^{\Lambda_p}.$$
 Here we regard $e^{\Lambda_p}$ only a symbol indexed by $p$.

 We define the operators $H_{i,l}$ $(0\leq i\leq n,\text{ }l\ne 0)$, 
 $e^{\alpha}$ $(\alpha\in Q)$, $\partial_{\alpha_i}$ and $z^{H_{i,0}}$ 
 $(0\leq i\leq n)$ on $V(p)_n$ as follows:\\
 for $v\otimes e^{\beta}e^{\Lambda_p}=H_{i_1,-k_1}\cdots H_{i_N,-k_N}v_0
 \otimes e^{\beta}e^{\Lambda_p}\in V(p)_n$ $(\beta=\sum_{k=0}^n
 m_k\alpha_k\in Q)$,
 $$H_{i,l}(v\otimes e^{\beta}e^{\Lambda_p})=(H_{i,l}v)\otimes
 e^{\beta}e^{\Lambda_p},$$
 $$e^{\alpha}(v\otimes e^{\beta}e^{\Lambda_p})=v\otimes
 (e^{\alpha}e^{\beta})e^{\Lambda_p},$$
 $${\partial_{\alpha_i}}(v\otimes e^{\beta}e^{\Lambda_p})=
 {\langle h_i,\beta+\Lambda_p\rangle}v\otimes e^{\beta}e^{\Lambda_p},$$
 $$z^{H_{i,0}}(v\otimes e^{\beta}e^{\Lambda_p})
  =z^{\langle h_i,\beta+\Lambda_p\rangle}
  \kappa^{\frac12 \sum_{k=0}^n\langle h_i,m_k\alpha_k\rangle m_{ik}} 
  v\otimes e^{\beta}e^{\Lambda_p}.$$
 $$
 d_1(v\otimes e^{\beta}e^{\Lambda_p})=(-\sum_{s=1}^{N}k_s
 -\frac{(\beta|\beta)}{2}
 -(\beta|\Lambda_p))v\otimes
 e^{\beta}e^{\Lambda_p},
 $$
 $$
 d_2(v\otimes e^{\beta}e^{\Lambda_p})=m_0(v\otimes e^{\beta}e^{\Lambda_p}).
 $$
 
The following lemma is easy.

\begin{lemma}
 As operators on $V(p)_n$,
 $$e^{\alpha_i}e^{\alpha_j}=(-1)^{\langle h_i,\alpha_j\rangle}
  e^{\alpha_j}e^{\alpha_i}$$
 $$q^{\partial_{\alpha_i}}e^{\alpha_j}=q^{\langle h_i,\alpha_j\rangle}
  e^{\alpha_j}q^{\partial_{\alpha_i}}$$
 $$z^{H_{i,0}}e^{\alpha_j}=z^{\langle h_i,\alpha_j\rangle}
  \kappa^{\frac12 \langle h_i,\alpha_j\rangle m_{ij}}e^{\alpha_j}
  z^{H_{i,0}}.$$
\end{lemma}

\begin{prop}
 Assume $c=1$ and $n>1$. Then
 for each $p$, the following action gives a $\tor$-module structure 
 on $V(p)_n$:
 $$q^{\frac12 c}\mapsto q^{\frac12},$$
 $$q^{d_1}\mapsto q^{d_1}$$
 $$q^{d_2}\mapsto q^{d_2}$$
 $$E_i(z)\mapsto\exp(\sum_{k\geq 1}\frac{H_{i,-k}}{[k]}(q^{-1/2}z)^k)
         \exp(\sum_{k\geq 1}-\frac{H_{i,k}}{[k]}(q^{1/2}z)^{-k})
         e^{\alpha_i}z^{H_{i,0}+1},$$
 $$F_i(z)\mapsto\exp(\sum_{k\geq 1}-\frac{H_{i,-k}}{[k]}(q^{1/2}z)^k)
         \exp(\sum_{k\geq 1}\frac{H_{i,k}}{[k]}(q^{-1/2}z)^{-k})
         e^{-\alpha_i}z^{-H_{i,0}+1},$$
 $$K_i^+(z)\mapsto \exp((q-q^{-1})\sum_{k\geq 1}
   H_{i,k}z^{-k})q^{\partial_{\alpha_i}},$$
 $$K_i^-(z)\mapsto \exp(-(q-q^{-1})\sum_{k\geq 1}
   H_{i,-k}z^{k})q^{-\partial_{\alpha_i}}$$
 for $0\leq i\leq n$.
\end{prop}

The proof will be given in Appendix.\\

It is easy to see the following lemma.
\begin{lemma}$V(p)_n$ is a cyclic $\tor$-module:
 $V(p)_n=\tor (v_0\otimes e^{\Lambda_p}).$
\end{lemma}

\begin{lemma}
 $V(p)_n$ has level $1$ as a $\affa$-module and as a $\affb$-module.
\end{lemma}

\begin{pf}
 It is clear that $V(p)_n$ is a level $1$ $\affa$-module. The center of $\affb$
 is $\prod_{k=0}^nK_i$. By the definition it acts as the scalar $q$ on 
 $V(p)_n$.
\end{pf}
\section{On $U_q({\frak {sl}}_{2,tor})$}
\subsection{}
 In this section we assume that ${\frak g}={\frak {sl}}_2$. 
 We shall try to find
 finitely many generators of $U_q({\frak {sl}}_{2,tor})$.\\

 Let 
 $$E_i=E_{i,0},\quad F_i=F_{i,0},\quad q^{\pm h_i}=K_i^{\pm},\text{ for }
   i=0,1,$$
 $$E_{-1}=F_{0,1}K_0^{-},\quad F_{-1}=K_0^+E_{0,-1},\quad q^{\pm h_{-1}}
   =q^{\pm c}K_0^{\mp}.$$

 \begin{prop}
  $U_q({\frak sl}_{2,tor})$ is generated by $E_i$, $F_i$, $q^{\pm h_i}$ 
  $(i=-1,0,1)$, $q^{\pm \frac12 c}$, $q^{\pm d_1}$, $q^{\pm d_2}$.
 \end{prop}
 
 \begin{pf}
 Let $\cal A$ be the subalgebra of $U_q({\frak sl}_{2,tor})$ generated by
$E_i$, $F_i$, $q^{\pm h_i}$ $(i=-1,0,1)$, $q^{\pm \frac12 c}$, 
$q^{\pm d_1}$, $q^{\pm d_2}$. By the
definition we have $E_{0,-1}=q^{-h_0}F_{-1}$ and $F_{0,1}=E_{-1}q^{h_0}$.
Since
\begin{align*}
 [E_{0,0},F_{0,1}] & = \frac{1}{q-q^{-1}}q^{-\frac12 c}K^+_{0,1}\\
                   & = q^{-\frac12 c}q^{h_0}H_{0,1}
\end{align*}
and
\begin{align*}
 [E_{0,-1},F_{0,0}] & = -\frac{1}{q-q^{-1}}q^{\frac12 c}K^-_{0,-1}\\
                   & = q^{\frac12 c}q^{-h_0}H_{0,-1}.
\end{align*}
we deduce $H_{0,1}$ and $H_{0,-1}\in {\cal A}$. We recall (2.2.7)
\begin{equation*}
 [H_{i,k},E_{j,l}]=\frac{1}{k}[k\langle h_i,\alpha_j\rangle]q^{-\frac12 |k|c}
 \kappa^{-km_{ij}}E_{j,k+l},
\end{equation*}
\begin{equation*}
 [H_{i,k},F_{j,l}]=-\frac{1}{k}[k\langle h_i,\alpha_j\rangle]q^{\frac12 |k|c}
 \kappa^{-km_{ij}}F_{j,k+l}.
\end{equation*}
By these formulas we have $E_{i,k}$, $F_{i,k}\in {\cal A}$ for $i=0,1$,
$k\in {\Bbb Z}$ inductively.

On the other hand we know
\begin{align*}
 [E_{0,1},F_{0,1}]& = \frac{1}{q-q^{-1}}q^{-\frac12 c}K^+_{0,2}\\
                  & = \frac12(q-q^{-1})H_{0,1}^2+H_{0,2}.
\end{align*}
Therefore we get $H_{0,2}\in {\cal A}$. Similarly we have
$H_{i,l}\in {\cal A}$ for any $i,l$.

This completes the proof.
\end{pf}

 \begin{lemma}
  The following relations hold in $U_q({\frak sl}_{2,tor})$:
  \begin{equation}
   [q^{\pm h_i},q^{\pm h_j}]=0,
  \end{equation}
  \begin{equation}
   [q^{\pm d_i},q^{\pm h_j}]=0,
  \end{equation}
  \begin{equation}
   q^{d_1}E_jq^{-d_1}=q^{\delta_{j,-1}}E_j,
  \end{equation}
  $$
   q^{d_1}F_jq^{-d_1}=q^{-\delta_{j,-1}}F_j,
  $$
  \begin{equation}
   q^{d_2}E_jq^{-d_2}=q^{-1+\delta_{j,1}}E_j,
  \end{equation}
  $$
   q^{d_2}F_jq^{-d_2}=q^{1-\delta_{j,1}}F_j,
  $$
  \begin{equation} 
   q^{h_i}E_jq^{-h_i}=q^{a_{ij}}E_j,
  \end{equation}
  $$
   q^{h_i}F_jq^{-h_i}=q^{-a_{ij}}F_j,
  $$
  where 
  \begin{equation*}
   (a_{ij})_{-1\leq i,j\leq 1}=
   \begin{pmatrix}
     2 & -2 &  2 \\
    -2 &  2 & -2 \\
     2 & -2 &  2
   \end{pmatrix},
  \end{equation*}
  \begin{equation}
   [E_i,F_j]=\delta_{ij}\frac{q^{h_i}-q^{-h_i}}{q-q^{-1}},
   \qquad\text{for $|i-j|\leq 1$,}
  \end{equation}
  \begin{equation}
   E_{-1}^3F_1-q^{-2}[3]E_{-1}^2F_1E_{-1}
   +q^{-4}[3]E_{-1}F_1E_{-1}^2-q^{-6}F_1E_{-1}^3=0,
  \end{equation}
  \begin{equation}
   E_{1}^3F_{-1}-q^{2}[3]E_{1}^2F_{-1}E_{1}
   +q^{4}[3]E_{1}F_{-1}E_{1}^2-q^{6}F_{-1}E_{1}^3=0,  
  \end{equation}
  \begin{equation}
   {}F_{-1}^3E_1-q^{-2}[3]F_{-1}^2E_1F_{-1}
   +q^{-4}[3]F_{-1}E_1F_{-1}^2-q^{-6}E_1F_{-1}^3=0,
  \end{equation}
  \begin{equation}
   {}F_{1}^3E_{-1}-q^{2}[3]F_{1}^2E_{-1}F_{1}
   +q^{4}[3]F_{1}E_{-1}F_{1}^2-q^{6}E_{-1}F_{1}^3=0,
  \end{equation}
  \begin{equation}
   E_i^3E_j-[3]E_i^2E_jE_i+[3]E_iE_jE_i^2-E_jE_i^3=0,
   \qquad\text{for $|i-j|=1$,}
  \end{equation}
  \begin{equation}
   {}F_i^3F_j-[3]F_i^2F_jF_i+[3]F_iF_jF_i^2-F_jF_i^3=0,
   \qquad\text{for $|i-j|=1$,}
  \end{equation}
  \begin{equation}
   E_{-1}E_1-q^2E_1E_{-1}=0,
  \end{equation}
  \begin{equation}
   {}F_{-1}F_1-q^2F_1F_{-1}=0.
  \end{equation}
 \end{lemma}

 \begin{pf}
  By the definition of $\tor$ and $E_i$, $F_i$ and $q^{h_i}$ it is easy to
  check these relations.
 \end{pf}

 Let ${\cal U}$ be an associative algebra over ${\Bbb Q}(q)$ generated by
 $E_i$, $F_i$, $q^{\pm h_i}$ $(i=-1,0,1)$, $q^{\pm \frac12 c}$, $q^{\pm d_1}$,
 $q^{\pm d_2}$ with relations
 (4.1.1)-(4.1.14). Then we have,

 \begin{cor}
  There is a canonical surjective algebra homomorphism $\Psi:{\cal U}\to 
  U_q({\frak {sl}}_{2,tor})$.
 \end{cor}

 \begin{rem}
  $\Psi$ has a highly nontrivial kernel. It is important to decide it.
  For example the following formulas holds in $U_q({\frak {sl}}_{2,tor})$:
  $$
  \kappa^{m_{01}}E_{0,0}E_{1,-1}-q^{-2}\kappa^{m_{01}}E_{1,-1}E_{0,0}
  =q^{-2}E_{0,-1}E_{1,0}-E_{1,0}E_{0,-1},
  $$  
  $$
  E_{0,-1}=q^{-h_0}F_{-1}
  $$
  and
  $$
  E_{1,-1}=\frac{\kappa^{-m_{01}}}{[-2]}[F_{-1}F_0-q^{-2}F_0F_{-1},E_1].
  $$
  Therefore we have
  \begin{align*}
   X & = \frac{1}{[-2]}E_0[F_{-1}F_0-q^{-2}F_0F_{-1},E_1]-
         \frac{q^{-2}}{[-2]}[F_{-1}F_0-q^{-2}F_0F_{-1},E_1]E_0\\
     & \qquad
         -q^{-2}q^{-h_0}F_{-1}E_1-E_1q^{-h_0}F_{-1}\\
     & = 0
  \end{align*}
  in $U_q({\frak sl}_{2,tor})$. Thus $X\in \text{Ker}\Psi$.
  But, as an element of $\cal U$, $X$ is not 
  equal to $0$. 
 \end{rem}
\subsection{}
Let 
$$E_{0^*}=F_{1,1}K_1^{-},\qquad F_{0^*}=K_1^+E_{1,-1},\qquad
  q^{\pm h_{0^*}}=q^{\pm c}K_1^{\mp}.$$
\begin{prop}
 The subalgebra generated by $E_i$, $F_i$ and $q^{h_i}$ for $i=0,1$,
 $E_{0^*}$, $F_{0^*}$, $q^{\pm h_{0^*}}$, $q^{\pm c}$, $q^{\pm d_1}$, 
 $q^{\pm d_2}$ is equal 
 to $U_q({\frak {sl}}_{2,tor})$. That is, they are generators of 
 $U_q({\frak {sl}}_{2,tor})$. Moreover these generators satisfy relations 
 similar to the ones in  Lemma 4.1.2.
\end{prop}

\begin{pf}
 This proposition is proved in the same way as Proposition 4.1.1 and Lemma
 4.1.2. 
\end{pf}

 We have immediately the following lemma.
 \begin{lemma}
 Let $U_q^{(1)}(\widehat{\frak {sl}}_2)$ be the subalgebras generated by $E_i$,
 $F_i$, $q^{h_i}$ for $i=1,0^*$ and $q^{\pm (d_1+d_2)}$, 
 $U_q^{(2)}(\widehat{\frak {sl}}_2)$ the 
 subalgebras generated by $E_i$, $F_i$,
 $q^{h_i}$ for $i=0,1$ and $q^{\pm (d_1+d_2)}$, 
 $U_q^{(3)}(\widehat{\frak {sl}}_2)$ 
 the subalgebras generated by
 $E_{i}$, $F_{i}$, $q^{\pm h_{i}}$ for $i=0^*,-1$ and $q^{\pm (d_1+d_2)}$,
 and $U_q^{(4)}(\widehat{\frak {sl}}_2)$ the 
 subalgebras generated by
 $E_i$, $F_i$, $q^{\pm h_i}$ for $i=0,-1$ and $q^{\pm (d_1+d_2)}$. 
 Then $U_q^{(i)}$ $(i=1,2,3,4)$ are isomorphic to 
 $U_q(\widehat{\frak {sl}}_2)$.
\end{lemma}

Those four algebras are schematically visualized by Fig. 1.\\

\begin{center}
\begin{picture}(0,0)%
\includegraphics{din.ps}%
\end{picture}%
\setlength{\unitlength}{0.01250000in}%
\begingroup\makeatletter\ifx\SetFigFont\undefined
\def\x#1#2#3#4#5#6#7\relax{\def\x{#1#2#3#4#5#6}}%
\expandafter\x\fmtname xxxxxx\relax \def\y{splain}%
\ifx\x\y   
\gdef\SetFigFont#1#2#3{%
  \ifnum #1<17\tiny\else \ifnum #1<20\small\else
  \ifnum #1<24\normalsize\else \ifnum #1<29\large\else
  \ifnum #1<34\Large\else \ifnum #1<41\LARGE\else
     \huge\fi\fi\fi\fi\fi\fi
  \csname #3\endcsname}%
\else
\gdef\SetFigFont#1#2#3{\begingroup
  \count@#1\relax \ifnum 25<\count@\count@25\fi
  \def\x{\endgroup\@setsize\SetFigFont{#2pt}}%
  \expandafter\x
    \csname \romannumeral\the\count@ pt\expandafter\endcsname
    \csname @\romannumeral\the\count@ pt\endcsname
  \csname #3\endcsname}%
\fi
\fi\endgroup
\begin{picture}(130,133)(280,598)
\put(290,605){\makebox(0,0)[lb]{\smash{\SetFigFont{12}{14.4}{rm}$-1$}}}
\put(405,710){\makebox(0,0)[lb]{\smash{\SetFigFont{12}{14.4}{rm}$1$}}}
\put(405,605){\makebox(0,0)[lb]{\smash{\SetFigFont{12}{14.4}{rm}$0^*$}}}
\put(410,655){\makebox(0,0)[lb]{\smash{\SetFigFont{12}{14.4}{rm}$U_q^{(1)}$}}}
\put(340,600){\makebox(0,0)[lb]{\smash{\SetFigFont{12}{14.4}{rm}$U_q^{(3)}$}}}
\put(280,655){\makebox(0,0)[lb]{\smash{\SetFigFont{12}{14.4}{rm}$U_q^{(4)}$}}}
\put(295,710){\makebox(0,0)[lb]{\smash{\SetFigFont{12}{14.4}{rm}$0$}}}
\put(340,720){\makebox(0,0)[lb]{\smash{\SetFigFont{12}{14.4}{rm}$U_q^{(2)}$}}}
\end{picture}

\end{center}
\begin{center}
{}Fig. 1.
\end{center}
Let $U_q(\frak {sl}_{2})_{(i)}$ $(i=-1,0,1,0^*)$ be the sub algebra of 
$U_q(\frak {sl}_{2,tor})$ generated by
$E_i$, $F_i$, $q^{\pm h_i}$. All 
$U_q(\frak {sl}_{2,(i)})$ $(i\in\{-1,0,1,0^*\})$ are isomorphic to
$U_q(\frak {sl}_2)$. The upper left
circle in Fig. 1 means $U_q(\frak {sl}_{2})_{(0)}$, the upper right one
means $U_q(\frak {sl}_{2})_{(1)}$, the lower left one means 
$U_q(\frak {sl}_{2})_{(-1)}$ and the lower right one means
$U_q(\frak {sl}_{2})_{(0^*)}$. The diagram 
$$\circ^i\Longleftrightarrow\circ^j\qquad (i,j\in\{-1,0,1,0^*\})$$
means the algebra generated by $U_q(\frak {sl}_{2})_{(i)}$ and
$U_q(\frak {sl}_{2})_{(j)}$ is isomorphic to $U_q(\widehat{\frak {sl}}_2)$.
{}For example $\circ^0\Longleftrightarrow\circ^1$ means the algebra 
generated by $U_q(\frak {sl}_{2})_{(0)}$ and $U_q(\frak {sl}_{2})_{(1)}$ which
we call $U_q^{(2)}(\widehat{\frak {sl}}_2)$ is isomorphic to
$U_q(\widehat{\frak {sl}}_2)$. The meaning of the diagram 
$$\circ^i==\Rightarrow\circ^j$$
is as follows: In the algebra generated by $U_q(\frak {sl}_{2})_{(i)}$ and
$U_q(\frak {sl}_{2})_{(j)}$, the following relations hold
\begin{equation}
q^{h_i}E_jq^{-h_i}=q^2E_j,\qquad q^{h_j}E_iq^{-h_j}=q^2E_i,
\end{equation}
\begin{equation}
q^{h_i}F_jq^{-h_i}=q^{-2}F_j,\qquad q^{h_j}F_iq^{-h_j}=q^{-2}F_i,
\end{equation}
\begin{equation}
E_i^3F_j-q^{-2}[3]E_i^2F_jE_i+q^{-4}[3]E_iF_jE_i^2-q^{-6}F_jE_i^3=0,
\end{equation}
$$
E_j^3F_i-q^{2}[3]E_j^2F_iE_j+q^{4}[3]E_jF_iE_j^2-q^{6}F_iE_j^3=0,
$$
\begin{equation}
{}F_i^3E_j-q^{-2}[3]F_i^2E_jF_i+q^{-4}[3]F_iE_jF_i^2-q^{-6}E_jF_i^3=0,
\end{equation}
$$
{}F_j^3E_i-q^{2}[3]F_j^2E_iF_j+q^{4}[3]F_jE_iF_j^2-q^{6}E_iF_j^3=0,
$$
\begin{equation}
E_iE_j-q^2E_jE_i=0,
\end{equation}
\begin{equation}
{}F_iF_j-q^2F_jF_i=0.
\end{equation}
\appendix
\section{ }
\subsection{Proof of Proposition 3.2.2 and 3.4.2.}
 {}For the proof we rewrite the defining relation of $\tor$ 
 generating function level.\\

 \begin{equation}
 q^{\pm\frac12 c}\text{ are central,}
\end{equation}\label{defrel1}
\begin{equation}
 K_i^+K_i^{-}=K_i^{-}K_i^+=1,
\end{equation}
\begin{equation}
 K_i^{\pm}(z)K_j^{\pm}(w)=K_j^{\pm}(w)K_i^{\pm}(z)
\end{equation}\label{rb}
\begin{equation}
 \theta_{-\langle h_i,\alpha_j\rangle}
 (q^{-c}\kappa^{m_{ij}}\frac{z}{w})K_i^-(z)K_j^+(w)=
 \theta_{-\langle h_i,\alpha_j\rangle}
 (q^{c}\kappa^{m_{ij}}\frac{z}{w})K_j^+(w)
 K_i^-(z)
\end{equation}\label{rc}
\begin{equation}
 q^{d_1}K_j^{\pm}(z)q^{-d_1}=K_j^{\pm}(q^{-1}z),
\end{equation}
\begin{equation}
 [q^{d_2},K_j^{\pm}(z)]=0,
\end{equation}
\begin{equation}
 q^{d_1}E_j(z)q^{-d_1}=E_j(q^{-1}z),
\end{equation}
$$
 q^{d_1}F_j(z)q^{-d_1}=F_j(q^{-1}z),
$$
\begin{equation}
 q^{d_2}E_j(z)q^{-d_2}=q^{\delta_{j0}}E_j(z),
\end{equation}
$$
 q^{d_2}F_j(z)q^{-d_2}=q^{-\delta_{j0}}F_j(z),
$$
\begin{equation}
 K_i^+(z)E_j(w)
 =\theta_{-\langle h_i,\alpha_j\rangle}
 (q^{-\frac12 c}\kappa^{-m_{ij}}\frac{w}{z})
 E_j(w)K_i^+(z)
 \end{equation}\label{rd}
$$
 K_i^-(z)E_j(w)
 =\theta_{\langle h_i,\alpha_j\rangle}
 (q^{-\frac12 c}\kappa^{m_{ij}}\frac{z}{w})
 E_j(w)K_i^-(z)
$$
$$
 K_i^+(z)F_j(w)
 =\theta_{\langle h_i,\alpha_j\rangle}
 (q^{\frac12 c}\kappa^{-m_{ij}}\frac{w}{z})
 {}F_j(w)K_i^+(z)
$$
\begin{equation*}
 K_i^-(z)F_j(w)
 =\theta_{-\langle h_i,\alpha_j\rangle}
 (q^{\frac12 c}\kappa^{m_{ij}}\frac{z}{w})
 {}F_j(w)K_i^-(z)
\end{equation*}
\begin{equation}
 [E_i(z),F_j(w)]=\delta_{i,j}\frac{1}{q-q^{-1}}
 \{\delta(q^c\frac{w}{z})K_i^+(q^{\frac12 c}w)-\delta(q^c\frac{z}{w})K_i^
 -(q^{\frac12 c}z)\}
\end{equation}\label{re}
\begin{equation}
 (\kappa^{m_{ij}}z-q^{\langle h_i,\alpha_j\rangle}w)E_i(z)E_j(w)
 =(q^{\langle h_i,\alpha_j\rangle}\kappa^{m_{ij}}z-w)E_j(w)E_i(z)
\end{equation}\label{rf}

$$
 (\kappa^{m_{ij}}z-q^{-\langle h_i,\alpha_j\rangle}w)F_i(z)F_j(w)
 =(q^{-\langle h_i,\alpha_j\rangle}\kappa^{m_{ij}}z-w)F_j(w)F_i(z)
$$
\begin{equation}
 \sum_{\sigma\in {\frak S}_m}\sum_{r=0}^m(-1)^r
 \begin{bmatrix}m\\r\end{bmatrix}_q
 E_i(z_{\sigma(1)})\cdots E_i(z_{\sigma(r)})E_j(w)E_i(z_{\sigma(r+1)})\cdots
 E_i(z_{\sigma(m)})=0
\end{equation}\label{rg}
\begin{equation*}
  \sum_{\sigma\in {\frak S}_m}\sum_{r=0}^m(-1)^r
 \begin{bmatrix}m\\r\end{bmatrix}_q
{} F_i(z_{\sigma(1)})\cdots F_i(z_{\sigma(r)})F_j(w)F_i(z_{\sigma(r+1)})\cdots
{} F_i(z_{\sigma(m)})
\end{equation*}

where $i\ne j$ and $m=1-\langle h_i,\alpha_j\rangle$.\\

In these formulas we denote $\theta_m(z)=
\frac{zq^m-1}{z-q^m}$ for
$m\in {\Bbb Z}$, $\delta(z)=\sum_{k\in {\Bbb Z}}z^k$.\\

If Proposition 3.4.2 holds, then, from Lemma 3.2.1 and 3.4.1, we
have Proposition 3.2.2 by putting $h_i\mapsto \overline{h_i}$, $\alpha_i
\mapsto\overline{\alpha}_i$, $\kappa\mapsto 1$ and $z^{H_{i,0}}\mapsto
z^{\overline{\partial}_i}$. Therefore it is enough to show Proposition 3.4.2.\\

The relations (A.1.1), (A.1.2) and (A.1.3) are trivial. (A.1.4)
is just the commutation relations of Heisenberg algebra $S_n$. Therefore, by 
the definition of $V(p)_n$, it is clear that (A.1.4) holds. The relations
(A.1.5), (A.1.6) immediately follows from the definition of $d_1$ and $d_2$.\\ 

Let us show (A.1.7) and (A.1.8). Take 
$v\otimes e^{\beta}e^{\Lambda_p}\in V(p)_n$ where $\beta=\sum_{k=0}^n
 m_k\alpha_k\in Q$. Then we have
\begin{align*}
& q^{d_1}e^{\pm \alpha_j}z^{\pm H_{j,0}+1}q^{-d_1}
(v\otimes e^{\beta}e^{\Lambda_p})\\
& = q^{\mp (\alpha_j|\beta+\Lambda_p)-1}
    z^{\pm \langle h_j,\beta+\Lambda_p\rangle +1}
    \kappa^{\frac12 \sum_{k=0}^n\langle h_j,m_k\alpha_k\rangle m_{jk}}
    v\otimes e^{\pm \alpha_j}e^{\beta}e^{\Lambda_p}\\
& = (q^{-1}z)^{\pm \langle h_j,\beta+\Lambda_p\rangle +1}
    \kappa^{\frac12 \sum_{k=0}^n\langle h_j,m_k\alpha_k\rangle m_{jk}}
    v\otimes e^{\pm \alpha_j}e^{\beta}e^{\Lambda_p}\\
& = e^{\pm \alpha_j}(q^{-1}z)^{\pm H_{j,0}+1}(v\otimes e^{\beta}e^{\Lambda_p}).
\end{align*}

Therefore we have
\begin{align*}
q^{d_1}E_j(z)q^{-d_1} 
& = q^{d_1}\exp(\sum_{k\geq 1}\frac{H_{j,-k}}{[k]}q^{-\frac12 k}z^k)
    \exp(-\sum_{k\geq 1}\frac{H_{j,k}}{[k]}q^{-\frac12 k}z^{-k})
    e^{\alpha_j}z^{H_{j,0}+1}q^{-d_1}\\
& = \exp(\sum_{k\geq 1}\frac{H_{j,-k}}{[k]}q^{-\frac12 k}(q^{-1}z)^k)
    \exp(-\sum_{k\geq 1}\frac{H_{j,k}}{[k]}q^{-\frac12 k}(q^{-1}z)^{-k})
    e^{\alpha_j}(q^{-1}z)^{H_{j,0}+1}\\
& = E_j(q^{-1}z).
\end{align*}

Similarly we have 
$q^{d_1}F_j(z)q^{-d_1}=F_j(q^{-1}z)$ .\\
 
It is clear that 
\begin{equation}
q^{d_2}e^{\pm \alpha_j}z^{\pm H_{j,0}+1}q^{-d_2}
(v\otimes e^{\beta}e^{\Lambda_p})
=q^{\pm \delta_{j0}}e^{\pm \alpha_j}z^{\pm H_{j,0}+1}
(v\otimes e^{\beta}e^{\Lambda_p}).
\end{equation}
 
{}From (A.1.13) and the fact that  $q^{d_2}$ commutes with $H_{j,k}$,
 we have (A.1.8).\\
  
We shall show (A.1.9).
We denote
 $$E_i^+(z)=\exp(\sum_{k\geq 1}\frac{H_{i,-k}}{[k]}q^{-\frac12 k}z^k),$$
 $$E_i^-(z)=\exp(-\sum_{k\geq 1}\frac{H_{i,k}}{[k]}q^{-\frac12 k}z^{-k}),$$
 $$F_i^+(z)=\exp(-\sum_{k\geq 1}\frac{H_{i,-k}}{[k]}q^{\frac12 k}z^k),$$
 $$F_i^-(z)=\exp(\sum_{k\geq 1}\frac{H_{i,k}}{[k]}q^{\frac12 k}w^{-k}).$$

Let us proof
\begin{equation*}
 K_i^+(z)E_j(w)
 =\theta_{-\langle h_i,\alpha_j\rangle}
 (q^{-\frac12}\kappa^{m_{-ij}}\frac{w}{z})
 E_j(w)K_i^+(z)
\end{equation*}

We have
\begin{align*}
 [(q-q^{-1})\sum_{k\ge 1}H_{i.k}z^{-k},\sum_{l\ge 1}\frac{H_{j.-l}}{[l]}
 q^{-\frac12 l}w^l]
 & = \sum_{k,l}(q-q^{-1})\frac{1}{[l]}[H_{i.k},H_{j.-l}]z^{-k}
     q^{-\frac12 l}w^l\\
 & = \sum_{k}(q-q^{-1})\frac{[k\langle h_i,\alpha_j\rangle]}{k}
     \kappa^{-km_{ij}}q^{-\frac12 k}(\frac{w}{z})^k\\
 & = \sum_{k}\frac{1}{k}(q^{k(\langle h_i,\alpha_j\rangle-\frac12)}
     -q^{k(-\langle h_i,\alpha_j\rangle-\frac12)})\kappa^{-km_{ij}}
     (\frac{w}{z})^k\\
 & = \text{log}\frac{1-q^{-\langle h_i,\alpha_j\rangle-\frac12}\kappa^{-m_{ij}}
     \frac{w}{z}}{1-q^{\langle h_i,\alpha_j\rangle-\frac12}\kappa^{-m_{ij}}
     \frac{w}{z}}
\end{align*}
and
$$
[(q-q^{-1})\sum_{k\ge 1}H_{i.k}z^{-k},-\sum_{l\ge 1}\frac{H_{j.l}}{[l]}
 q^{-\frac12 l}w^{-l}]=0.
$$
We recall Campbell-Hausdorff formula: let $A$ and $B$
be noncommutative operators and $C=[A,B]$. If $[C,A]=[C,B]=0$ then we have
$e^Ae^B=e^Ce^Be^A$.\\
 
By Campbell-Hausdorff formula we get
\begin{align*}
 & \exp((q-q^{-1})\sum_{k\ge 1}H_{i.k}z^{-k})E_j^+(w)E_j^-(w)\\
 & =\frac{1-q^{-\langle h_i,\alpha_j\rangle-\frac12}\kappa^{-m_{ij}}
     \frac{w}{z}}{1-q^{\langle h_i,\alpha_j\rangle-\frac12}\kappa^{-m_{ij}}
     \frac{w}{z}}
 E_j^+(w)E_j^-(w)\exp((q-q^{-1})\sum_{k\ge 1}H_{i.k}z^{-k}).
\end{align*}
On the other hand, by Lemma 3.4.1 we have
$$
 q^{\partial_{\alpha_i}}e^{\alpha_j}w^{H_{j,0}}
 =q^{\langle h_i,\alpha_j\rangle}e^{\alpha_j}w^{H_{j,0}}
 q^{\partial_{\alpha_i}}.
$$
Thus we get
\begin{align*}
 & K_i^+(z)E_j(w) \\
 & = \exp((q-q^{-1})\sum_{k\ge 1}H_{i.k}z^{-k})E_j^+(w)E_j^-(w)
     q^{\partial_{\alpha_i}}e^{\alpha_j}w^{H_{j,0}+1}\\
 & = \frac{1-q^{-\langle h_i,\alpha_j\rangle-\frac12}\kappa^{-m_{ij}}
     \frac{w}{z}}{q^{-\langle h_i,\alpha_j\rangle}-q^{-\frac12}\kappa^{-m_{ij}}
     \frac{w}{z}}E_j(w)K_i^+(z)\\
 & = \theta_{-\langle h_i,\alpha_j\rangle}
     (q^{-\frac12}\kappa^{-m_{ij}}\frac{w}{z})
     E_j(w)K_i^+(z).
\end{align*}
 
The other formulas in (A.1.9) can be checked by similar arguments.\\

Let us show (A.1.10). We have
 \begin{align*}
  & [-\sum_{k\geq 1}\frac{H_{i,k}}{[k]}q^{-\frac12 k}z^{-k},
  -\sum_{k\geq 1}\frac{H_{j,-k}}{[k]}q^{\frac12 k}w^k]\\
  & = \sum_{k\geq 1}\frac{1}{k[k]}[k\langle h_i,\alpha_j\rangle]
      \kappa^{-m_{ij}}(\frac{w}{z})^k\\
  & = \begin{cases}
       \text{log}\frac{1}{(1-q\frac{w}{z})
       (1-q^{-1}\frac{w}{z})},
          & \langle h_i,\alpha_j\rangle=2\quad (i=j),\\
       \text{log}(1-\kappa^{-m_{ij}}\frac{w}{z}),
          & \langle h_i,\alpha_j\rangle=-1,\\
       0,
          & \langle h_i,\alpha_j\rangle=0.
      \end{cases}
 \end{align*}
  
{}For example we will show in the case of $\langle h_i,\alpha_j\rangle=-1$.
By Campbell-Hausdorff formula we get
$$E_i^{\pm}(z)F_j^{\pm}(w)=F_j^{\pm}(w)E_i^{\pm}(z)$$
and
$$E_i^-(z)F_j^+(w)=(1-\kappa^{-m_{ij}}\frac{w}{z})F_j^+(w)E_i^-(z).$$

 On the other hand, by Lemma 3.4.1, we have
 $$z^{H_{i,0}}e^{-\alpha_j}
  =z\kappa^{\frac12 m_{ij}}e^{-\alpha_j}z^{H_{i,0}}.$$

 Therefore we get
 \begin{align*}
  E_i(z)F_j(w)
  & = E_i^+(z)E_i^-(z)
      e^{\alpha_i}z^{H_{i,0}+1}
      F_j^+(w)F_j^-(w)
      e^{-\alpha_j}w^{-H_{j,0}+1}\\
  & = (z\kappa^{\frac12 m_{ij}}-w\kappa^{-\frac12 m_{ij}})
      E_i^+(z)F_j^+(w)E_i^-(z)F_j^-(w)\\
  & \qquad \times
      e^{\alpha_i}e^{-\alpha_j}
      z^{H_{i,0}+1}
      w^{-H_{j,0}+1}.
 \end{align*}

 By a similar argument we have
 \begin{align*}
{}  F_j(w)E_i(z)
  & = -(w\kappa^{-\frac12 m_{ij}}-z\kappa^{\frac12 m_{ij}})
      E_i^+(z)F_j^+(w)E_i^-(z)F_j^-(w)\\
  & \qquad \times
      e^{\alpha_i}e^{-\alpha_j}
      z^{H_{i,0}+1}
      w^{-H_{j,0}+1}.
 \end{align*}

Therefore we get
 $$[E_i(z),F_j(w)]=0.$$

Similarly one can check the other formulas.\\

We will show (A.1.11). We have
\begin{align*}
  & [-\sum_{k\geq 1}\frac{H_{i,k}}{[k]}q^{-\frac12 k}z^{-k},
  \sum_{k\geq 1}\frac{H_{j,-k}}{[k]}q^{-\frac12 k}w^k]\\
  & = \sum_{k\geq 1}\frac{1}{k[k]}[k\langle h_i,\alpha_j\rangle]
      \kappa^{-m_{ij}}(\frac{w}{z})^k\\
  & = \begin{cases}
       \text{log}(1-\frac{w}{z})
       (1-q^{-2}\frac{w}{z}),&
       i=j,\\
       \text{log}\frac{1}{1-q^{-1}\kappa^{-m_{ij}}\frac{w}{z}},&
       \langle h_i,\alpha_j\rangle=-1,\\
       0,&
       \langle h_i,\alpha_j\rangle=0.
      \end{cases}
 \end{align*}

{}For example let us show in the case of $\langle h_i,\alpha_j\rangle=-1$.
If $\langle h_i,\alpha_j\rangle\ne-1$ one can show the formula by a similar
argument.
By Lemma 3.4.1 we get
\begin{align*}
& (z\kappa^{m_{ij}}-q^{\langle h_i,\alpha_j\rangle}w)E_i(z)E_j(w)\\
& = (z\kappa^{m_{ij}}-q^{-1}w)E_i^+(z)E_i^-(z)
    e^{\alpha_i}z^{H_{i,0}+1}
    E_j^+(w)E_j^-(w)
    e^{\alpha_j}w^{H_{j,0}+1}\\
& = \frac{z\kappa^{m_{ij}}-q^{-1}w}{z\kappa^{\frac12 m_{ij}}
    (1-q^{-1}\kappa^{-m_{ij}}\frac{w}{z})}
    E_i^+(z)E_j^+(w)E_i^-(z)E_j^-(w)\\
& \qquad \times
    e^{\alpha_i}e^{\alpha_j}
    z^{H_{i,0}+1}w^{H_{j,0}+1}\\
& = \kappa^{\frac12 m_{ij}}E_i^+(z)E_j^+(w)E_i^-(z)E_j^-(w)
    e^{\alpha_i}e^{\alpha_j}
    z^{H_{i,0}+1}w^{H_{j,0}+1}.
\end{align*}
On the other hand we have
\begin{align*}
& (q^{\langle h_i,\alpha_j\rangle}\kappa^{m_{ij}}z-w)E_j(w)E_i(z)\\
& = \frac{q^{-1}\kappa^{m_{ij}}z-w}
    {w\kappa^{-\frac12 m_{ij}}
     (1-q^{-1}\kappa^{m_{ij}}\frac{z}{w})}
    E_j^+(w)E_i^+(z)E_j^-(w)E_i^-(z)\\
& \qquad \times
    e^{\alpha_j}e^{\alpha_i}
    w^{H_{j,0}+1}z^{H_{i,0}+1}\\
& = \kappa^{\frac12 m_{ij}}E_i^+(z)E_j^+(w)E_i^-(z)E_j^-(w)
    e^{\alpha_i}e^{\alpha_j}
    z^{H_{i,0}+1}w^{H_{j,0}+1}.
\end{align*}
Thus we have $(z\kappa^{m_{ij}}-q^{\langle h_i,\alpha_j\rangle}w)E_i(z)E_j(w)
=(q^{\langle h_i,\alpha_j\rangle}\kappa^{m_{ij}}z-w)E_j(w)E_i(z)$.

The formula $ (z\kappa^{m_{ij}}-q^{-\langle h_i,\alpha_j\rangle}w)F_i(z)F_j(w)
 =(q^{-\langle h_i,\alpha_j\rangle}\kappa^{m_{ij}}z-w)F_j(w)F_i(z)$ 
is proved similarly.\\

Let us prove (A.1.12). Assume that $\langle h_i,\alpha_j\rangle=-1$.
This is the most complicated case. The other cases can be proved similarly.\\

We have the following formulas:
\begin{align*}
 & E_i(z_1)E_i(z_2)E_j(w)\\
 & = \frac{(z_1-z_2)(z_1-q^{-2}z_2)\kappa^{-m_{ij}}}
          {z_1z_2(1-\frac{q^{-1}\kappa^{-m_{ij}}w}{z_1})
           (1-\frac{q^{-1}\kappa^{-m_{ij}}w}{z_2})}
 E_i^+(z_1)E_i^+(z_2)E_j^+(w)E_i^-(z_1)E_i^-(z_2)E_j^-(w)\\
 & \qquad \times
 e^{2\alpha_i}e^{\alpha_j}
 z_1^{H_{i,0}+1}z_2^{H_{i,0}+1}w^{H_{j,0}+1},
\end{align*}
\begin{align*}
 & E_i(z_1)E_j(w)E_i(z_2)\\
 & = \frac{(z_1-z_2)(z_1-q^{-2}z_2)}
          {z_1w(1-\frac{q^{-1}\kappa^{-m_{ij}}w}{z_1})
           (1-\frac{q^{-1}\kappa^{m_{ij}}z_2}{w})}
 E_i^+(z_1)E_i^+(z_2)E_j^+(w)E_i^-(z_1)E_i^-(z_2)E_j^-(w)\\
 & \qquad \times
 e^{2\alpha_i}e^{\alpha_j}
 z_1^{H_{i,0}+1}z_2^{H_{i,0}+1}w^{H_{j,0}+1},
\end{align*}
\begin{align*}
 & E_j(w)E_i(z_1)E_i(z_2)\\
 & = \frac{(z_1-z_2)(z_1-q^{-2}z_2)\kappa^{m_{ij}}}
          {w^2(1-\frac{q^{-1}\kappa^{m_{ij}}z_1}{w})
           (1-\frac{q^{-1}\kappa^{m_{ij}}z_2}{w})}
 E_i^+(z_1)E_i^+(z_2)E_j^+(w)E_i^-(z_1)E_i^-(z_2)E_j^-(w)\\
 & \qquad \times
 e^{2\alpha_i}e^{\alpha_j}
 z_1^{H_{i,0}+1}z_2^{H_{i,0}+1}w^{H_{j,0}+1}.
\end{align*}
Therefore it is enough to show that
\begin{equation}
\begin{split}
 & \sum_{\sigma\in{\frak S}_2}
 (z_{\sigma(1)}-z_{\sigma(2)})(z_{\sigma(1)}-q^{-2}z_{\sigma(2)})\\
 & \qquad\times\{
 \frac{\kappa^{-m_{ij}}}
    {z_{\sigma(1)}z_{\sigma(2)}(1-\frac{q^{-1}\kappa^{-m_{ij}}w}
     {z_{\sigma(1)}})
     (1-\frac{q^{-1}\kappa^{-m_{ij}}w}{z_{\sigma(2)}})}
 -\frac{q+q^{-1}}
    {z_{\sigma(1)}w(1-\frac{q^{-1}\kappa^{-m_{ij}}w}{z_{\sigma(1)}})
     (1-\frac{q^{-1}\kappa^{m_{ij}}z_{\sigma(2)}}{w})}\\
 & \qquad
 +\frac{\kappa^{m_{ij}}}
    {w^2(1-\frac{q^{-1}\kappa^{m_{ij}}z_{\sigma(1)}}{w})
     (1-\frac{q^{-1}\kappa^{m_{ij}}z_{\sigma(2)}}{w})}
  \}\\
 & = 0.
\end{split}
\end{equation}
This identity is proved by a direct calculation.\\
Thus the proposition is proved.

\begin{ack}
The author would like to thank Masaki Kashiwara, Tetsuji Miwa and Michio Jimbo
for their encouragement and valuable discussions. He also thanks Kenji Iohara,
Norio Suzuki, Eric Vasserot and Hiroshi Yamada for stimulating discussions.
\end{ack}


\begin{thebibliography}{[MFT50]}
\bibitem[Be]{}J.~Beck,
{\it Braid group action and quantum affine algebras},
Commun. Math. Phys. 165 (1994), 555-568.
\bibitem[B]{}S.~Berman,
{\it On generators and relations for certain involutory subalgebras Kac-Moody
Lie algebras},
Commun. Alg. 17(12) (1989), 3165-3185.
\bibitem[BM]{}S.~Berman and R.~Moody,
{\it Lie algebras graded by finite root systems and the intersection matrix
algebras of Slodowy},
Inv. Math. (1992), 323-347.
\bibitem[C]{}I.~Cherednik,
{\it Double-affine Hecke algebras, Knizhink-Zamolodchikov equations and
Macdonald's operators},
Intern. Math. Research Notices 9 (1992), 171-180.
\bibitem[CP]{}V.~Chari and A.~Pressly,
{\it Quantum affine algebras and affine Hecke algebras},
Pacific Jour. of Math. 174 No.2 (1996), 295-326.
\bibitem[D1]{}V.~Drinfeld,
{\it Quantum groups},
Proc. ICM Berkeley 1 (1986), 789-820.
\bibitem[D2]{}V.~Drinfeld,
{\it A new realization of Yangians and quantized affine algebras},
Sov. Math. Dokl. 36 (1988), 212-216.
\bibitem[EM]{}S.~Eswara Rao and R.~V.~Moody,
{\it Vertex representations for N-toroidal Lie algebras and a Generalization
of the Virasoro algebra},
Commun. Math. Phys. 159 (1994), 239-264.
\bibitem[FJ]{}I.~Frenkel and N.~Jing,
{\it Vertex representations of quantum affine algebras},
Proc. Natl. Acad. Sci. USA 85 (1988), 9373-9377.
\bibitem[FK]{}I.~Frenkel and V.~Kac,
{\it Basic Representations of affine Lie algebras and dual Resonance models},
Inv. Math. 62 (1980), 23-66.
\bibitem[GKV]{}V.~Ginzburg, M.~Kapranov and E.~Vasserot,
{\it Langlands reciprocity for algebraic surface},
Math. Res. Lrtt. 2 (1995), 147-160.
\bibitem[K]{}Y.~Koyama,
{\it Staggered polarization of vertex model with 
$U_q(\widehat{\frak sl}(n))$-symmetry},
Commun. Math. Phys. 164 (1994), 277-291.
\bibitem[MEY]{}R.~V.~Moody, S.~Eswara Rao and T~.Yokomuma,
{\it Toroidal Lie algebras and vertex representations},
Geom. Dedicata 35 (1990), 283-307.
\bibitem[S]{}P.~Slodowy,
{\it Beyond Kac-Moody algebras, and inside},
Canadian Math. Soc. Conference Proceedings Vol.5 (1986), 361-371.
\bibitem[V]{}E.~Vasserot,
{\it Private communication.}
\bibitem[VV]{}M.~Varagnolo and E.~Vasserot,
{\it Schur duality in the toroidal setting}
Commun. Math. Phys. 182 (1996), 469-484.
\bibitem[Y]{}H.~Yamada,
{\it Extended affine Lie algebra and their vertex representations},
Publ. RIMS. 25 (1989), 587-603. 
\end{thebibliography}
\end{document}